\documentclass[11pt]{article}
\usepackage{tikz}
\usepackage[T1]{fontenc}
\usepackage{amsmath,amssymb,amsthm,mathrsfs}
\usepackage{microtype}
\usepackage{hyperref}
\allowdisplaybreaks
\newtheorem{theorem}{Theorem}[section]
\newtheorem{lemma}[theorem]{Lemma}

\newtheorem{claim}{Claim}

\newtheorem{conjecture}[theorem]{Conjecture}

\newcommand{\eps}{\varepsilon}
\newcommand{\one}{\mathbf 1}

\newcommand{\lam}{\lambda_{\max}}

\begin{document}
\title{Maxima of the Q-index: Forbidden an odd wheel
}
\author{Hongyu Wang$^a$,  \quad Yue Ma$^b$\\
\small $^{a}$ School of  Sciences,\\
\small Southwest Petroleum University, 
\small Chengdu, Sichuan 610500, China.\\
\small $^{b}$School of Mathematics and Statistics,\\
\small Nanjing University of Science and Technology,
\small Nanjing, Jiangsu 210094, China.\\
\small $^a$ {why71@mail.ustc.edu.cn} $^b$yma@njust.edu.cn}

\date{}

\maketitle

\begin{abstract}
In 2013, Freitas, Nikiforov and Patuzzi (EJLA, 2013) proposed the problem of determining the maximum signless Laplacian spectral radius of an $n$-vertex $F$-free graph. Recently, Yuan, Zhai and Zhang (LAA, 2026) determined the maximum signless Laplacian spectral radius of graphs without wheels of sufficiently large order. In this paper, we consider the problem of forbidding a single wheel $W_{2k+1}=K_1\vee C_{2k}$. We show that, for every fixed integer $k\ge 3$ and sufficiently large $n$, the unique extremal graph maximizing the signless Laplacian spectral radius among all $n$-vertex $W_{2k+1}$-free graphs is $K_k\vee(K_2\cup(n-k-2)K_1)$.
\end{abstract}

\section{Introduction}
All graphs considered in this paper are finite and simple, and subgraph containment is not necessarily induced.
For a graph $G=(V(G),E(G))$, we call $|V(G)|$ the \emph{order} of $G$ and $e(G):=|E(G)|$ the \emph{size} of $G$.
Let $A(G)$ be its adjacency matrix, $D(G)$ the diagonal matrix of vertex degrees, and $Q(G)=D(G)+A(G)$ its signless Laplacian. The largest eigenvalue of $Q(G)$ is called the $Q$-\emph{index} of $G$, denoted by $q(G)$.
For a fixed graph $F$, the signless Laplacian spectral extremal problem asks for the maximum of $q(G)$ among all $n$-vertex $F$-free graphs. This problem was proposed by de Freitas, Nikiforov and Patuzzi~\cite{FreitasNikiforovPatuzzi2013}.

For integers $s,t$, let $K_s$, $C_s$ and $P_s$ be the clique, cycle and path of order $s$, respectively. Let $K_{s,t}=sK_1\vee tK_1$ be the complete bipartite graph.
Let $W_{s+1}:=K_1\vee C_s$ be the wheel of order $s+1$. We call a wheel an odd wheel if its order is odd. 
For $1\le k<n$, let $S_{n,k}:=K_k\vee(n-k)K_1$ be the complete split graph of order $n$, and for $k\le n-2$ let $S_{n,k}^{+}$ be obtained from $S_{n,k}$ by adding one edge to its independent set. 
To be more specific, $S_{n,k}^{+}:=K_k\vee\bigl(K_2\cup(n-k-2)K_1\bigr)$.

Several results are known for the maximum $Q$-index of a graph with forbidden subgraphs,
including cliques~\cite{HeJinZhang2013},
paths~\cite{NikiforovYuan2014},
cycles~\cite{FreitasNikiforovPatuzzi2013, Yuan2014, NikiforovYuan2015}, 
$K_{2,t}$'s~\cite{FreitasNikiforovPatuzzi2016},
books~\cite{ChenJinZhang2025}, fans~\cite{WangZhai2023},
linear forests~\cite{ChenLiuZhang2020,ChenLiuZhang2023},
vertex-disjoint cliques~\cite{ChenJinZhang2026, ZhangWang2024, ZhaoYouZeng2026},
intersecting odd cycles~\cite{ZhaoHuangGuo2021, ChenLiuZhang2024}
and so on.

For wheels, Zhao, Huang and Lin~\cite{ZhaoHuangLin2021} proved that $S_{n,2}$ is the unique $n$-vertex wheel-free graph with maximum $Q$-index for $n\ge4$. Chen and Ma~\cite{ChenMa2023} determined the maximum
$Q$-index of graphs without both $W_5$ and $C_6$. Recently, Yuan, Zhai and Zhang~\cite{YuanZhaiZhang2026} determined the maximum $Q$-index of graphs without all wheels of sufficiently large order. For a fixed odd wheel, the adjacency spectral extremal problem was studied by Cioab\u{a}, Desai and Tait~\cite{CioabaDesaiTait2022} and Zhang~\cite{Zhang2025}. A related adjacency spectral problem with a fixed number of edges was considered in~\cite{ChenGaoLi2026}.

For general forbidden graphs, Zheng, Li and Su~\cite{ZhengLiSu2026} proved a signless Laplacian spectral version of the Erd\H{o}s--Stone--Simonovits theorem for graphs with chromatic number at least three. Moreover, Zheng, Li and Li~\cite{ZhengLiLi2026} proved an exact signless Laplacian spectral Tur\'an theorem for color-critical graphs with chromatic number at least four. 
In particular, their result determines the maximum $Q$-index of graphs without a fixed even wheel. However, it does not apply to any odd wheel $W_{2k+1}$.

In this paper, we determine the maximum $Q$-index of $n$-vertex graphs containing no $W_{2k+1}$ for fixed $k\ge3$ and sufficiently large $n$. 

\begin{theorem}\label{thm:main}
For every fixed integer $k\ge3$, there is a constant $N(k)$ such that, for any $n\ge N(k)$, every $n$-vertex $W_{2k+1}$-free graph $G$ satisfies $q(G)\le q(S_{n,k}^{+}),$ with equality if and only if $G\cong S_{n,k}^{+}$.
\end{theorem}

The rest of the paper is arranged as follows. In Section~2, we present some preliminary results, including structural lemmas and a Schur-complement estimate. We prove the extremal result for graphs with sublinearly many large degree vertices in Section~3. Section~4 establishes a spectral stability theorem. Finally, Section~5 completes the proof of Theorem~\ref{thm:main}.

\section{Preliminaries}
\subsection{Basic notations and estimates}
Given a graph $G$, a vertex $v\in V(G)$  and sets $X,Y\subseteq V(G)$. 
Let $K[X]$ be the complete graph on $X$ and $K[X,Y]$ be the complete bipartite graph with partition set $X$ and $Y$ if they are disjoint. Let $G[X]=G\cap K[X]$ and $G[X,Y]=G\cap K[X,Y]$.
Denote by $e(X)$ the number of edges of $G[X]$ and by $e(X,Y)$ the number of edges of $G[X,Y]$. 
Denote by $N_G(v)$ the neighborhood of $v$ in $G$ (we may write $N(v)$ for short) and $d_G(v):=|N_G(v)|$ the degree of $v$ in $G$. Also, $N_X(v):=N_G(v)\cap X$, $d_X(v):=|N_X(v)|$, $N_X(Y):=\bigcup_{v\in Y}N_X(v)$ and $\Delta(G):=\max_{v\in G}d_{G}(v)$. 

\begin{lemma}\label{lem:benchmark}
Let $k\ge3$ and $n\ge4k$. The graph $S_{n,k}^{+}$ is $W_{2k+1}$-free. If $q_F=q(S_{n,k}^{+})$ and $q_S=q(S_{n,k})$, then
$$
n<n+2(k-1)-\frac{2k(k-1)}{n}<q_S<q_F<n+2k.
$$
Moreover, $q_F$ is the largest root of
\begin{align*}
P_{n,k}(x)&=x^3-(n+3k)x^2+(kn+4k^2+2n-4)x-2k^2(k+1).
\end{align*}
\end{lemma}
\begin{proof}
A vertex outside the $K_k$-part of $S_{n,k}^{+}$ has degree at most $k+1<2k$ and cannot be a hub (i.e. the only vertex of degree $2k$ in $W_{2k+1}$). For a hub in $K_k$, a $2k$-cycle in its neighborhood would use at least $k+1$ vertices outside the $K_k$-part and thus at least two edges outside. This is a contradiction since there is only one edge outside the $K_k$-part in $S_{n,k}^{+}$. 
The equitable signless Laplacian quotient of $S_{n,k}^{+}$, with the classes corresponding to the $K_k$-part, the two outside edge endpoints, and the remaining vertices, is
$$
\begin{pmatrix}n+k-2&2&n-k-2\\k&k+2&0\\k&0&k\end{pmatrix}.
$$
Its positive quotient Perron eigenvector lifts to a positive eigenvector of $Q(S_{n,k}^{+})$, yielding the displayed characteristic polynomial. 
The analogous quotient of $S_{n,k}$ has characteristic polynomial $p(x)=x^2-(n+2(k-1))x+2k(k-1)$. As adding an outside edge strictly increases the $Q$-index of the connected graph $S_{n,k}$, we have $q_F>q_S$.
One can check that $p(n)<0$, $P_{n,k}(n+2k)>0$ and $P_{n,k}$ is increasing on $[n,+\infty)$. Hence $q_S>n$ and $q_F<n+2k$. Together with $p(q_S)=0$, we get $n+2(k-1)-q_S=2k(k-1)/q_S<2k(k-1)/n$. The proof is done.
\end{proof}

\subsection{Structural lemmas}

\begin{lemma}\label{lem:star}
Let $k\ge 3$ and $G$ be a $W_{2k+1}$-free graph. Pick $u\in V(G)$, $S\subseteq N(u)$ with $|S|=k-1$, and let $T\subseteq V(G)\backslash(S\cup\{u\})$ be a vertex set whose vertices are adjacent to every vertex of $S\cup\{u\}$. then:
\begin{enumerate}
\item[(1)] If  $|T|\ge k+1$, then $e(T)\le1$.
\item[(2)] if $|T|\ge k$, then $|N(z)\cap T|\le1$ holds for any $z\in N(u)\backslash (S\cup T)$.
\end{enumerate}
\end{lemma}
\begin{proof}
For (1), suppose to the contrary that $|T|\ge k+1$ and $e(T)\ge2$.
Clearly, we have either $P_3\cup (k-2)K_1\subseteq G[T]$ or $2K_2\cup (k-3)K_1\subseteq G[T]$.
We obtain $k-1$ vertex-disjoint paths covering $k+1$ vertices in $T$ in both cases. 
Hence, we can use the $k-1$ vertices in $S$ to connect the endpoints of these paths cyclically, forming a $C_{2k}$ in $G[N(u)]$. This is a contradiction since $G$ is $W_{2k+1}$-free.

For (2), suppose $|T|\ge k$ and $|N(z)\cap T|\ge2$ for some $z\in N(u)\backslash (S\cup T)$.
Clearly, we have $P_3\cup (k-2)K_1\subseteq G[T\cup\{z\}]$ and $z$ is not an endpoint but the middle vertex of the $P_3$. Similarly as the proof of (1), we then get a $C_{2k}$ in $G[N(u)]$, a contradiction.
\end{proof}

\begin{lemma}\label{lem:common}
Let $G$ be a $W_{2k+1}$-free graph and $A,B$ be disjoint vertex sets in $G$ with $|A|\ge k$. Put $M=\sum_{a\in A}|B\setminus N(a)|$. Then
$$
e(B)\le\left(\frac{k-1}{2}+\frac{kM}{|A|}\right)|B|.
$$
\end{lemma}
\begin{proof}
Let $T=B\cap\bigcap_{s\in S}N(s)$.
By the pigeonhole principle, we can pick a set $S\subseteq A$ with $|S|=k$ and $|B\setminus T|\le\sum_{s\in S}|B\setminus N(s)|\le \frac{k}{|A|}M$. If some $v\in T$ has $k$ neighbors in $G[T]$, these $k$ vertices alternating the $k$ vertices in $S$ form a $C_{2k}$ in $N(v)$, a contradiction. Thus $\Delta(G[T])\le k-1$. Therefore,
$$
e(B)\le e(T)+|B\setminus T|\,|B|\le\left(\frac{k-1}{2}+\frac{kM}{|A|}\right)|B|.
$$
\end{proof}

For an $n$-vertex graph $G$ and a number $x\le n-1$, Let $L_G(x)=\{v:d_G(v)\ge x\}$.
Fix $0<\eps<1/2$ and $\alpha=\frac{\eps}{2k}$. Define
$$
L=L_G(\eps n), \quad R=V(G)\setminus L, \quad U=L_G(n-1-\alpha n), \quad W=L\setminus U.
$$
\begin{lemma}\label{lem:nearu}
There exists $N_0=N_0(k,\eps)$ such that for every
$n\ge N_0$ and every $n$-vertex $W_{2k+1}$-free graph $G$,
then:
\begin{enumerate}
\item[(1)] $d_U(v)\le k-1$ for every vertex $v\in L$.
\item[(2)] If $d_U(u)= k-1$ for some $u\in U$, then $d_W(u)=0$.
\item[(3)] For every vertex $u\in U$, $d_L(u)\le D_0:=\left\lceil (k-1)(4/\eps)^k\right\rceil.$
\end{enumerate}
\end{lemma}
\begin{proof}
For (1), suppose that $v\in L$ is adjacent to $k$ vertices $u_1,\ldots,u_k\in U$. Since
$d_G(v)\ge\eps n$ and each $u_j$ has at most $n-d_G(u_j)\le\alpha n+1$ nonneighbors (including itself),
$$
\left|N(v)\cap\bigcap_{j=1}^kN(u_j)\right|
\ge \eps n-k(\alpha n+1)>k
$$
for $n\ge N_0(k,\eps)$. Hence, 
the vertices $u_1,\ldots,u_k$ together with $k$ vertices in $N(v)\cap\bigcap_{j=1}^kN(u_j)$ form an alternating
$C_{2k}$ inside $N(v)$, a contradiction.

For (2), let $S\subseteq N_U(u)$ with $|S|=k-1$, and suppose that $z\in N_W(u)$.
Note that $z\in W\subseteq L$ and $\{u\}\cup S\subseteq U$, similarly as the previous proof,
$$
\left|N(z)\cap N(u)\cap\bigcap_{s\in S}N(s)\right|
\ge\eps n-k(\alpha n+1)>k
$$
for $n\ge N_0(k,\eps)$.
Hence, the $k$ vertices of $S\cup\{z\}$ and $k$ vertices in $N(z)\cap N(u)\cap\bigcap_{s\in S}N(s)$
form an alternating $C_{2k}$ inside $N(u)$, a contradiction, too.

For (3), every vertex $v\in N_L(u)=N(u)\cap L$ satisfies
$$
d_{N(u)}(v)
\ge d_G(v)-|V(G)\setminus N(u)|
\ge\eps n-(n-d_G(u))
\ge\eps n-\alpha n-1
\ge\frac{\eps n}{2}.
$$
for $n\ge N_0(k,\eps)$.
Since $G[N(u)]$ is $C_{2k}$-free, 
it is $K_{k,k}$-free. Thus, every $k$-subset of $N(u)$ can be contained in the neighborhoods of at most $k-1$ vertices in $N_L(u)$. This gives
$$
d_L(u)\binom{\lfloor\eps n/2\rfloor}{k}
\le \sum_{v\in N_L(u)}\binom{d_{N(u)}(v)}{k}
\le(k-1)\binom{d_G(u)}{k}.
$$
Note that
$\lfloor\eps n/2\rfloor-k+1\ge\eps n/4$ for $n\ge N_0(k,\eps)$.
Then
$$
\binom{\lfloor\eps n/2\rfloor}{k}
\ge\frac{(\eps n/4)^k}{k!},
\qquad\binom{d_G(u)}{k}\le\frac{n^k}{k!}.
$$
It follows that
$d_L(u)\le(k-1)(4/\eps)^k\le D_0$.
\end{proof}

We may also use the classical result of Bondy and Simonovits~\cite{BondySimonovits1974} that, for every fixed integer $k\ge 3$, every $C_{2k}$-free graph on $m$ vertices has at most $O_k(m^{1+1/k})=o_k(m^2)$ edges.

\subsection{A Schur-complement estimate}
For integers $s,t>0$, let $J_{s\times t}$ be the $s\times t$ matrix whose entries are all $1$'s. Write $J_t:=J_{t\times t}$ and $\one_t=J_{t\times 1}$ for short. Also, $I_t$ is the $t\times t$ unit matrix. We use $\lam(M)$ to denote the spectral radius of a matrix $M$.

In this subsection, suppose $|L|=o(n)$, $q=q(G)>n$, and let $x=(x_v)$ be a nonnegative $q$-eigenvector normalized by $\max_vx_v=1$. Order vertices by $L\sqcup R$. We write $x^{\mathsf T}=(x_L^{\mathsf T},x_R^{\mathsf T})$ and write 
$$
Q(G)=\begin{pmatrix}A&B\\B^{\mathsf T}&C\end{pmatrix},\qquad
P=B(qI_{|R|}-C)^{-1}B^{\mathsf T}.
$$
Put $a_0=(1-2\eps)^{-1}$ and $r=16a_0e(R)/n^2$. 

\begin{lemma}\label{lem:schur0}
For all sufficiently large $n$, $P\le (1+r)J_{|L|}$ entrywise.
\end{lemma}
\begin{proof}
We write $Q_R=Q(G[R])$ for short.
Since $0\le C\le |L|I_{|R|}+Q_R$ entrywise and both matrices have spectral radius smaller than $q$, the comparison of their nonnegative Neumann series and simple calculations give
\begin{align*}
(qI_{|R|}-C)^{-1}&\le((q-|L|)I_{|R|}-Q_R)^{-1}\\
&=\frac {1}{q-|L|}(I_{|R|}+Q_R((q-|L|)I_{|R|}-Q_R)^{-1})
\end{align*}
entrywise. Hence
\begin{equation}\label{eq:schurP}
\begin{aligned}
P &=B(qI_{|R|}-C)^{-1}B^{\mathsf T}\\
&\le \frac {1}{q-|L|}B(I_{|R|}+Q_R((q-|L|)I_{|R|}-Q_R)^{-1})B^{\mathsf T}\\
&=\frac{1}{q-|L|}BB^{\mathsf T}+\frac{1}{q-|L|}BQ_R((q-|L|)I_{|R|}-Q_R)^{-1}B^{\mathsf T}\\
&\le J_{|L|}+\frac{1}{q-|L|}BQ_R((q-|L|)I_{|R|}-Q_R)^{-1}B^{\mathsf T}
\end{aligned}
\end{equation}
entrywise. The last inequality holds since the entries of $\frac{1}{q-|L|}BB^{\mathsf T}$ are at most $\frac{n-|L|}{q-|L|}<1$. 

We next estimate the second term. 
Let $b=(b_v)$ be any column of $B^{\mathsf T}$. Then $0\le b\le\one$. Set $z=(z_v)=((q-|L|)I_{|R|}-Q_R)^{-1}b$.
$((q-|L|)I_{|R|}-Q_R)^{-1}$ is entrywise nonnegative by its Neumann series, so $z\ge 0$ entrywise.
Choose $u\in R$ such that $z_u=\max_{v\in R}z_{v}$.
Note that $b_u\le1$ and $z_v\le z_u$ for every $v\in R$.
Consider the row of $u$ in $((q-|L|)I_{|R|}-Q_R)z=b$, we get
$$
(q-|L|)z_u=b_u+d_{R}(u)z_u+\sum_{v\in N_R(u)}z_v\le1+2d_{R}(u)z_u\le1+2\Delta(G[R])z_u.
$$
This implies $z_u \le\frac1{q-|L|-2\Delta(G[R])}$. In other words, $z\le \frac1{q-|L|-2\Delta(G[R])}\one$.

This holds for every column of $B^{\mathsf T}$, which means
$$
((q-|L|)I_{|R|}-Q_R)^{-1}B^{\mathsf T}\le\frac{1}{q-|L|-2\Delta(G[R])} J_{|R|\times|L|}.
$$
Note that $0\le B\le J_{|L|\times |R|}$ entrywise, and $J_{|L|\times |R|}Q_RJ_{|R|\times |L|}=4e(R)J_{|L|}$. It holds
\begin{align*}
BQ_R((q-|L|)I_{|R|}-Q_R)^{-1}B^{\mathsf T}&\le \frac{1}{q-|L|-2\Delta(G[R])}J_{|L|\times|R|}Q_R J_{|R|\times|L|}\\
&=\frac{4e(R)}{q-|L|-2\Delta(G[R])}J_{|L|}.
\end{align*}
This inequation, together with \eqref{eq:schurP}, gives $P\le (1+\frac{4e(R)}{(q-|L|)(q-|L|-2\Delta(G[R]))})J_{|L|}$. It remains to prove $\frac{4e(R)}{(q-|L|)(q-|L|-2\Delta(G[R]))}<r=16a_0e(R)/n^2$.
This is trivial since $q-|L|>\frac{n}{2}$ and $q-|L|-2\Delta(G[R])>\frac{1-2\eps}{2}n=\frac{n}{2a_0}$ for large $n$.
\end{proof}

\begin{lemma}\label{lem:schur}
Let $X=\sum_{v\in L}x_v$, $\delta_R(v)=|R\setminus N_G(v)|,$ and
$$
\mathcal M=A(G[L])+J_{|L|}-nI_{|L|}+\operatorname{diag}(d_G(v):v\in L).
$$
We have
\begin{align}
(q-n)x_L&\le\mathcal M x_L+rX\one,\label{eq:schurpoint}\\
(q-n)X&\le\sum_{v\in L}(2d_{L}(v)-\delta_R(v))x_v+r|L|X,\label{eq:schursum}\\
q-n&\le\lam(\mathcal M)+r|L|.\label{eq:schurlambda}
\end{align}
Furthermore, for a vertex $u$ with $x_u=1$, we have $u\in L$ and
\begin{equation}\label{eq:schurmax}
q-n\le d_{L}(u)-\delta_R(u)-|L|+\sum_{v\in N_{L}(u)}x_v+(1+r)X.
\end{equation}
\end{lemma}
\begin{proof}
Since $qI-Q(G)$ is positive semidefinite, its Schur complement $qI-A-P$ is positive semidefinite and singular. Hence $q=\lam(A+P)$. Also, $(A+P)x_L=qx_L$. 
Note that $\mathcal M=A(G[L])+J-nI+\operatorname{diag}(d_G(v):v\in L)=A+J-nI$. We then have $\mathcal Mx_L=(q-n)x_L+X\one-Px_L$. By Lemma~\ref{lem:schur0}, we soon get \eqref{eq:schurpoint}.
It can be checked that the row sum at $v\in L$ of $\mathcal M$ equals $2d_{L}(v)-\delta_R(v)$. By \eqref{eq:schurpoint}, this gives \eqref{eq:schursum}. 

By Lemma~\ref{lem:schur0}, $A+P\le A+(1+r)J_{|L|}$ entrywise. Both matrices are nonnegative and symmetric, so Perron--Frobenius monotonicity gives $q\le\lam(A+(1+r)J_{|L|})$. Therefore $q-n\le\lam(\mathcal M+rJ_{|L|})\le\lam(\mathcal M)+\lam(rJ_{|L|})=\lam(\mathcal M)+r|L|$, which completes the proof of \eqref{eq:schurlambda}.
If $u\in R$, its eigenvalue equation would imply $q\le2d_G(u)<n$, which is a contradiction. Hence $u\in L$. Consider $u\in L$ with $x_u=1$ in \eqref{eq:schurpoint}, we get \eqref{eq:schurmax}.
\end{proof}

\section{Proof of the extremal case}
In this section, we deal with the case when the $n$-vertex $W_{2k+1}$-free graph contains only $o(n)$ vertices of degree at least $\Theta(n)$.
\begin{theorem}\label{thm:sublinear}
Fix $k\ge 3$ and $0<\eps<1/2$. If $G$ is an $n$-vertex $W_{2k+1}$-free graph with $q(G)\ge q(S_{n,k}^{+})$ and $|L_G(\eps n)|=o_{k,\eps}(n)$, then $G\cong S_{n,k}^{+}$ holds for all sufficiently large $n$.
\end{theorem}

We prove the theorem by several claims.
Throughout the proof,
we consider an $n$-vertex $W_{2k+1}$-free graph $G$ with $n\ge N$, $q\ge q(S_{n,k}^{+})$ and $|L|=o(n)$,
where we write $q=q(G)$ and use $L,R,U,W$ defined before. 
Moreover, define $\delta(v)=n-1-d_G(v)$ and $\delta_R(v)=|R\setminus N(v)|$ for $v\in L$.
Choose a nonnegative $q$-eigenvector $x$ with $\max_{v\in V(G)}x_v=1$, and put $X=\sum_{v\in L}x_v$.
By Lemma~\ref{lem:benchmark},
\begin{equation}\label{eq:delta-lower-revised}
q-n>2(k-1)-\frac{2k(k-1)}{n}>0.
\end{equation}
Lemma~\ref{lem:schur} ensures that a vertex $u$ with $x_u=1$ belongs to $L$.
Recall $\alpha=\eps/(2k)$, $a_0=(1-2\eps)^{-1}$ and $r=16a_0e(R)/n^2$, and set
$$
D_0=\left\lceil (k-1)(4/\eps)^k\right\rceil,
\qquad b_0 =\frac{1+a_0}{\alpha}, \qquad b_1=b_0\frac{|L|}{n}=o(1).
$$

\begin{claim}\label{c:mass-local1}
$u\in U$.
\end{claim}

\begin{proof}
For large $n$, $R$ is nonempty as $|L|=o(n)$.
Pick a vertex $v_m\in R$ with $x_{v_m}=\max_{v\in R}x_{v}$.
The eigenvalue equation at $v_m$, together with $x_v\le1$ on $L$, gives
$$
(q-d_G(v_m))x_{v_m}
 =\sum_{v\in N_L(v_m)}x_v+\sum_{v\in N_R(v_m)}x_v
 \le |L|+d_{R}(v_m)x_{v_m}.
$$
Since $d_{R}(v_m)\le d_G(v_m)<\eps n$ and $q>n$, it follows that
$$
x_{v_m}\le \frac{|L|}{q-d_G(v_m)-d_{R}(v_m)}
 \le \frac{|L|}{(1-2\eps)n}=a_0\frac{|L|}{n}.
$$
Consequently,
\begin{equation}\label{eq:mass0}
\sum_{v\in V(G)}x_v
 \le |L|+|R|x_{v_m}\le(1+a_0)|L|.
\end{equation}
As $x_u=1$, the eigenvalue equation at $u$ implies
$$
q-d_G(u)=\sum_{v\in N(u)}x_v
 \le\sum_{v\in V(G)}x_v.
$$
Note that $\delta(u)=q-d_G(u)-(1+q-n)<q-d_G(u)$. Thus
\begin{equation}\label{eq:bu-revised}
0\le \delta(u)<(1+a_0)|L|=o(n).
\end{equation}
In particular, for large $n$, we have $\delta(u)\le\alpha n$, so $u\in U$.
\end{proof}

\begin{claim}\label{c:mass-local2}
$x_v\le b_1$ for $v\in W$.
\end{claim}
\begin{proof}
By Lemma~\ref{lem:nearu}, since $n$ is sufficiently large,
\begin{equation}\label{eq:dhu-revised}
d_L(u)\le D_0.
\end{equation}
For $v\in W$, we have $\delta(v)>\alpha n$. 
Applying the eigenvalue equation at $v\in W$ and using \eqref{eq:mass0} and $q-d_G(v)=1+q-n+\delta(v)> \delta(v)>\alpha n$, we get
\begin{equation*}
0\le x_v=\frac{\sum_{i\in N(v)}x_i}{q-d_G(v)}\le\frac{(1+a_0)|L|}{\alpha n}
 =b_0\frac{|L|}{n}=b_1=o(1)
\end{equation*}
holds for every $v\in W$.
\end{proof}
\begin{claim}\label{c:mass-local3}
$e(R)=o(n^2)$. In particular, $r=o(1)$.
\end{claim}
\begin{proof}
Every $W_{2k+1}$-free graph has $C_{2k}$-free neighborhoods. In particular,
$G[R\cap N(u)]$ is $C_{2k}$-free and has $o(n^2)$ edges.
The set $R\setminus N(u)$ has $\delta_R(u)\le \delta(u)=o(n)$ vertices,
whose degrees are all less than $\eps n$. Hence
$$
e(R)\le e(G[R\cap N(u)])+\eps n \delta_R(u)=o(n^2).
$$
It follows that $r=16a_0e(R)/n^2=o(1)$. 
\end{proof}

\begin{claim}\label{lem:linear-bootstrap}
$r|L|=o(1)$.
\end{claim}

\begin{proof}
By Claim~\ref{c:mass-local3}, we can assume $|L|\to\infty$.
Let $X=\sum_{v\in L}x_v$. 
The fact that $\delta_R(u)+q-n\ge0$ and $\sum_{v\in N_L(u)}x_v\le d_L(u)$, together with \eqref{eq:schurmax} and \eqref{eq:dhu-revised}, implies
\begin{equation}\label{eq:mass}
|L|-X\le 2D_0+rX\le 2D_0+r|L|.
\end{equation}
By Claim~\ref{c:mass-local3} and \eqref{eq:mass}, we have $X=|L|-o(|L|)$. Put
$$
A_* =\{v\in L:x_v\ge\tfrac12\},\qquad X_W=\sum_{i\in W}x_i.
$$
The inequality $X\le |A_*|+\frac12(|L|-|A_*|)$ gives $|A_*|\ge2X-|L|\ge |L|/2$ for large $n$, so in particular $|A_*|\ge k$. By \eqref{eq:mass0}, for each $v\in A_*$, the eigenvalue equation gives
$$
\delta(v) x_v\le(q-d_G(v))x_v
 =\sum_{i\in N(v)}x_{i}\le(1+a_0)|L|.
$$
Since $x_v\ge1/2$, we obtain $\delta(v)\le2(1+a_0)|L|=o(n)$. Thus $A_*\subseteq U$ for all sufficiently large $n$.

By Claim~\ref{c:mass-local2},
$$
(1-b_1)|W|\le\sum_{v\in W}(1-x_v)
 \le |L|-X\le2D_0+r|L|.
$$
Taking $n$ sufficiently large that $b_1\le1/2$, we have
\begin{equation}\label{eq:wcount-revised}
|W|\le4D_0+2r|L|
 =4D_0+\frac{32a_0|L|}{n^2}e(R).
\end{equation}

Define
$$
M=\sum_{v\in L}\delta_R(v)x_v.
$$
By Lemma~\ref{lem:nearu} (1), every vertex has at most $k-1$ neighbors in $U$. Therefore
$e(U,W)=\sum_{v\in W}d_U(v)\le (k-1)|W|$. Since $L=U\cup W$,
\begin{align}
\sum_{v\in L}d_L(v)x_v
&=\sum_{v\in U}(d_U(v)+d_W(v))x_v
  +\sum_{v\in W}(d_U(v)+d_W(v))x_v\notag\\
&\le (k-1)\sum_{v\in U}x_v+e(U,W)
  +(k-1)\sum_{v\in W}x_v+|W|\sum_{v\in W}x_v\notag\\
&\le (k-1)X+(k-1)|W|+|W|X_W.\label{eq:M-revised0}
\end{align}
Here $d_W(v)\le |W|$ for $v\in W$ and $x_v\le1$ for $v\in U$.
Using \eqref{eq:schursum},
\eqref{eq:delta-lower-revised} and \eqref{eq:M-revised0}, we obtain
\begin{align}
M
&\le 2\sum_{v\in L}d_L(v)x_v-(q-n) X+r|L|X\notag\\
&\le \frac{2k(k-1)}{n}X+2(k-1)|W|+2|W|X_W+r|L|X\notag\\
&\le \frac{2k(k-1)|L|}{n}+2(k-1)|W|+2b_1 |W|^2
   +\frac{16a_0|L|^2}{n^2}e(R).\label{eq:M-revised}
\end{align}

We apply Lemma~\ref{lem:common} with $A=A_*$ and $B=R$.
The sets are disjoint and $|A_*|\ge k$. Since $x_v\ge1/2$
for $v\in A_*$,
$$
\sum_{v\in A_*}\delta_R(v)\le2\sum_{v\in A_*}\delta_R(v)x_v\le2M.
$$
Thus, using $|A_*|\ge |L|/2$, $|R|\le n$ and \eqref{eq:M-revised},
\begin{align}
e(R)
&\le\left(\frac{k-1}{2}+\frac{k}{|A_*|}\sum_{v\in A_*}\delta_R(v)\right)|R|\notag\\
&\le\frac{(k-1)n}{2}+\frac{4kn}{|L|}M\notag\\
&\le\frac{(k-1)n}{2}+8k^2(k-1)
     +\frac{8kn|W|}{|L|}(k-1+b_1 |W|)
     +\frac {64ka_0|L|}{n} e(R).
\end{align}
Note that
$\frac{n}{|L|}b_1|W|^2=b_0|W|^2\le b_0|L||W|$. We conclude that
\begin{equation}\label{eq:bootstrap-short-revised}
\begin{aligned}
e(R)&\le\frac{(k-1)n}{2}+8k^2(k-1)\\
     &+8k\left(\frac{(k-1)n}{|L|}+b_0|L|\right)|W|
     +\frac {64ka_0|L|}{n} e(R).
\end{aligned}
\end{equation}
By \eqref{eq:wcount-revised} and
\eqref{eq:bootstrap-short-revised},
\begin{equation}\label{eq:bootstrap-final-revised}
\begin{aligned}
(1-\theta_n)e(R)
 &\le\frac{(k-1)n}{2}+8k^2(k-1)
 +32kD_0\left(\frac{(k-1)n}{|L|}+b_0|L|\right),\\
\text{where }& \theta_n
 =64ka_0\frac {|L|}{n}
 +256ka_0\left(\frac{|L|}{n}+b_0\frac{|L|^2}{n^2}\right).
\end{aligned}
\end{equation}
The constants  $k,\eps,\alpha, D_0, a_0, b_0$ are fixed.
$\frac{|L|}{n}\to0$ implies $\theta_n\to0$.
For large $n$, it holds $\theta_n\le1/2$. Since $1\le |L|\le n$,
the right side of \eqref{eq:bootstrap-final-revised}
is $O_{k,\eps}(n)$. Therefore, $e(R)=O_{k,\eps}(n)$, which implies
\begin{equation*}
r|L|=\frac{16a_0|L|}{n^2}e(R)=O_{k,\eps}(|L|/n)=o(1).
\end{equation*}
\end{proof}

\medskip

\begin{claim}\label{lem:integer-rigidity}
$\delta_R(u)=0$ and $U=L$. In particular, $d_L(u)=\Delta(G[L])=k-1$ and $|L|\ge k$.
\end{claim}

\begin{proof}
Write $d_U=d_U(u)$, $d_W=d_W(u)$  for short.
Since $u\in U$, Lemma~\ref{lem:nearu} (1)
gives $d_U\le k-1$. In \eqref{eq:schurmax}, the neighbors of $u$
in $U$ contribute at most $d_U$ to $\sum_{v\in N_L(u)}x_v$,
and those in $W$ contribute at most $b_1 d_W$ by Claim~\ref{c:mass-local2}.
Moreover, $|L|-X\ge(1-b_1)|W|$. 
Hence by \eqref{eq:schurmax},
\begin{equation}\label{eq:integer}
\begin{aligned}
q-n\le&2d_U+(1+b_1)d_W-(1-b_1)|W|-\delta_R(u)+r|L|,\\
 &d_U\le k-1,\quad d_W\le D_0,\quad |W|\ge d_W.
\end{aligned}
\end{equation}
If $d_U\le k-2$, then
$$
q-n\le2(k-2)+2b_1 d_W+r|L|
\le2k-4+2b_1 D_0+r|L|=2k-4+o(1),
$$
contrary to \eqref{eq:delta-lower-revised} for large enough $n$.
Therefore $d_U=k-1$. By Lemma~\ref{lem:nearu} (2), the fact $u\in U$ and
$d_U(u)=k-1$ force $d_W=0$. By \eqref{eq:integer} and \eqref{eq:delta-lower-revised},
$$
(1-b_1)|W|+\delta_R(u)\le2(k-1)-(q-n)+r|L|
 <\frac{2k(k-1)}{n}+r|L|=o(1).
$$
Here $|W|$ and $\delta_R(u)$ are nonnegative integers and $1-b_1\ge1/2$ for large $n$.
It follows that $|W|=\delta_R(u)=0$. This leads to $U=L$. Thus by Lemma~\ref{lem:nearu} (1),
\begin{equation*}
d_L(u)=\Delta(G[L])=k-1,\qquad |L|\ge k.
\end{equation*}
\end{proof}

\begin{claim}\label{lem:global-rigidity}
$e(R)\le1$ and $d_L(v)=k-1$, $\delta_R(v)=0$ for every $v\in L$.
\end{claim}
\begin{proof}
Since $U=L$, Lemma~\ref{lem:nearu} (1) gives $d_L(v)\le k-1$ for every $v\in L$. Define
$$
\sigma_v=2(k-1-d_L(v))+\delta_R(v)\ge 0,
\quad D=\sum_{v\in L}\sigma_v,
\quad \rho=2(k-1)-\lambda_{\max}(\mathcal M).
$$

The row sum of $\mathcal M$ at $v$ is
$2d_L(v)-\delta_R(v)=2(k-1)-\sigma_v$. 
It is easy to check by the definition of  $\mathcal M$ that $\mathcal M+nI_{|L|}$ is entrywise positive. 
So by the maximum-row-sum
inequality,
$\lambda_{\max}(\mathcal M)\le2(k-1)$, and then $\rho\ge0$.
By \eqref{eq:schurlambda}, \eqref{eq:delta-lower-revised}
and $r|L|=o(1)$,
\begin{equation}\label{eq:eta-small-revised}
0\le\rho\le2(k-1)-(q-n)+r|L|
 <\frac{2k(k-1)}{n}+r|L|=o(1).
\end{equation}

Let $y>0$ be an eigenvector of $\mathcal M$ for
$\lambda_{\max}(\mathcal M)=2(k-1)-\rho$, normalized by
$\max_{v\in L}y_v=1$. Pick $w\in L$ with $y_w=1$.
Recall that the row sum $\sum_{v\in L}\mathcal M_{wv}=2(k-1)-\sigma_w$. 
The eigenvector equation of the row of $w$ gives
$$
\rho=\sigma_w+\sum_{v\in L\setminus\{w\}}
     \mathcal M_{wv}(1-y_v).
$$
Every summand on the right is nonnegative and
$\mathcal M_{wv}\ge1$, so $y_v\ge1-\rho$ for every $v$,
which holds for $v=w$ as well.
On the other hand, the symmetry of $\mathcal M$ and its row-sum formula give
$$
(2(k-1)-\rho)\sum_{v\in L} y_v
 =\sum_{v\in L} (2(k-1)-\sigma_v)y_v,
\quad\text{i.e.}\quad
\rho\sum_{v\in L} y_v=\sum_{v\in L}\sigma_v y_v.
$$
Since $\rho=o(1)$, it holds $1-\rho\ge1/2$. Therefore
\begin{equation}\label{eq:defect-bound-revised}
\rho |L|\ge\rho\sum_{v\in L} y_v
 =\sum_{v\in L}\sigma_v y_v\ge(1-\rho)D\ge D/2,
\qquad\text{i.e.}\quad \rho\ge\frac{D}{2|L|}.
\end{equation}
Combining \eqref{eq:defect-bound-revised} with
\eqref{eq:eta-small-revised}, we get $D=o(|L|)$.

Let $B=\{v\in L:\sigma_v>0\}$. Because every positive $\sigma_v$ is an integer,
$|B|\le D$. By Claim~\ref{lem:integer-rigidity}, $\Delta(G[L])=k-1$, so
at most $|B|+(k-1)|B|\le kD=o(|L|)$ vertices lie in $B\cup N_L(B)$. 

Since $n$ is sufficiently large, we can pick a vertex
$v_0\in L\setminus(B\cup N_L(B))$.
Thus $\sigma_{v_0}=0$ and $\sigma_v=0$ for every $v\in N_L(v_0)$.
By the definition of $\sigma_v$, this implies $|N_L(v_0)|=d_L(v_0)=k-1$.
It also implies $\delta_R({v_0})=0$ and $\delta_R(v)=0$ for $v\in N_L(v_0)$, which means each vertex in $\{v_0\}\cup N_L(v_0)$ is adjacent to every vertex of $R$. 
Note that $N_L(v_0)\subseteq N(v_0)$, $R\subseteq V(G)\setminus(\{v_0\}\cup N_{L}(v_0))$, and $|R|=n-|L|=n-o(n)\ge k+1$ for large $n$.
By Lemma~\ref{lem:star} (1) (with $(S,T,u)$ in the lemma be $(N_L(v_0),R,v_0)$), we have $e(R)\le1$.

Now suppose $D\ge1$. By \eqref{eq:defect-bound-revised},
$\rho\ge\frac{1}{2|L|}$. Using \eqref{eq:schurlambda}
and $e(R)\le1$,
$$
q-n\le2(k-1)-\rho+r|L|
\le2(k-1)-\frac1{2|L|}+\frac{16a_0|L|}{n^2}.
$$
Together with \eqref{eq:delta-lower-revised}, this implies
$$
\frac12< 2k(k-1)\frac {|L|}n+16a_0\frac{|L|^2}{n^2},
$$
which is impossible for $|L|=o(n)$. 
Note that every $\sigma_v$ is a nonnegative integer, and so is $D$. Thus $D=0$, which immediately implies $\sigma_v=0$ for every $v\in L$. 
Therefore, $d_L(v)=k-1$ and $\delta_R(v)=0$ for every $v\in L$.
\end{proof}
\medskip

\begin{proof}[Proof of Theorem~\ref{thm:sublinear}]
Write $\ell=|L|$. By Claim~\ref{lem:integer-rigidity}, $\ell\ge k$.
By Claim~\ref{lem:global-rigidity}, $e(R)\le 1$.
Let $G^+=G$ if $e(R)=1$ and $G^+$ is obtained from $R$ by adding an edge in $R$
if $e(R)=0$. Hence
$$
G^+=G[L]\vee\bigl(K_2\cup(n-\ell-2)K_1\bigr)
$$
with $q(G)\le q(G^+)$. 
As $G[L]$ is $(k-1)$-regular by Claim~\ref{lem:global-rigidity}, the partition of $V(G)$ into $L$,
the two endpoints of the edge, and the remaining vertices is
equitable for $Q(G^+)$, with quotient
$$
\mathcal B_{\ell}=
\begin{pmatrix}
n-\ell+2k-2&2&n-\ell-2\\
\ell&\ell+2&0\\
\ell&0&\ell
\end{pmatrix}.
$$
The largest eigenvalue of $\mathcal B_{\ell}$ is $q(G^+)$.
Moreover, $G^+$ contains the connected spanning graph
$K_{\ell,n-\ell}$, which has $Q$-index $n$. Strict edge monotonicity therefore
gives $q(G^+)>n$. 
Let $x=q(G^+)>n>\ell+2$, we have $\det(xI-\mathcal B_\ell)=0$, thus
$$
f_\ell(x):=x^3-(n+\ell+2k)x^2+((\ell+2)n+4(k-1)(\ell+1))x-2\ell((\ell+2)k-\ell)=0.
$$
For sufficiently large $n$, one can check that $f_\ell$ is increasing on $[n,\infty)$ and $f_\ell(n)<0$. Thus its unique root on $[n,\infty)$ is exactly $x=q(G^+)$.

For $\ell=k$, $G[L]$ is necessarily $K_k$ and $G^+=S_{n,k}^{+}$.
Hence $f_k(q_F)=0$, where $q_F=q(S_{n,k}^{+})$.
By simple calculations,
$$
f_\ell(q_F)-f_k(q_F)
=(\ell-k)[(n+4k-4-q_F)q_F-2((k-1)\ell+k(k+1))].
$$
Since $n<q_F<n+2k$ by Lemma~\ref{lem:benchmark}, we have $(n+4k-4-q_F)q_F>(2k-4)n$ and $2((k-1)\ell+k(k+1))=o(n)$.
If $\ell>k$, then $f_\ell(q_F)=f_\ell(q_F)-f_k(q_F)>0$.
Since $f_\ell$ is increasing on $[n,\infty)$, it holds $q(G)\le q(G^+)<q_F,$ which is a contradiction.
Therefore, $\ell=k$ and $G[L]=K_k$.
Finally, if $e(R)=0$ then $G=S_{n,k}$,
whose $Q$-index is strictly smaller than $q_F$
by Lemma~\ref{lem:benchmark}.
Thus, $e(R)=1$ and $G\cong S_{n,k}^{+}$. The proof is done.
\end{proof}

\section{A spectral stability theorem}
For any graph $H$, we call a partition of its vertex set $V(H)=A\sqcup B$ a \emph{maximum cut} of $H$, if $e(A,B)=\max_{V(H)=A'\sqcup B'}e(A',B')$.
Let $V_1\sqcup V_2$ be a cut of some graph $H$, then for $i=1,2$ and some vertex $v\in V_i$, every vertex in $N_{V_{i}}(v)$ is called an \emph{internal neighbor} of $v$ and every vertex in $N_{V_{3-i}}(v)$ is called a \emph{crossing neighbor} of $v$.
We call edges in $E(H[V_1])\cup E(H[V_2])$ \emph{internal edges} and edges in $E(H[V_1,V_2])$ \emph{crossing edges}. We define internal nonneighbors, crossing nonneighbors, internal nonedges and crossing nonedges similarly.

\begin{lemma}\label{lem:cut}
Fix $k\ge 3$ and $0<\tau<1/2$. Let $G$ be a $W_{2k+1}$-free graph on $n$ vertices with $q(G)\ge n$ and $\Delta(G)\le(1-\tau)n$.
Then there is a maximum cut $A\sqcup B$ of $G$ with $|A|,|B|\ge\tau n/2$ and
$$
\begin{aligned}
\quad |A||B|-e(A,B)=O_{k,\tau}(n^{1+\frac{1}{k}}),
\quad e(A)+e(B)=O_{k,\tau}(n^{1+\frac{1}{k}}).
\end{aligned}
$$
Moreover, $d_{B}(v_a)\ge d_{A}(v_a)$ and $d_{A}(v_b)\ge d_{B}(v_b)$ for every vertex $v_a\in A$ and $v_b\in B$.
\end{lemma}
\begin{proof}
Without loss of generality, we only consider $G$ without isolated vertices.
By the maximum-row-sum inequality of $D(G)^{-1}Q(G)D(G)$ (see also~\cite{Merris1998, FengGu2009}), there is a vertex $u\in V(G)$ with
$$
q(G)\le d_{G}(u)+\frac{1}{d_{G}(u)}\sum_{v\in N(u)}d_G(v).
$$
Set $A=N(u)$, $B=V(G)\setminus A$, $a=|A|=d_{G}(u)$ and $b=|B|$. Since the right side of the inequation above is at most $a+\Delta(G)$, we have $a\ge q(G)-\Delta(G)\ge\tau n\ge k$. Also, $b\ge n-\Delta(G)\ge\tau n$. 
Apparently,
$$
n\le q(G)\le a+\frac1a\sum_{v\in A}d(v)=n+\frac{2e(A)+e(A,B)-ab}{a}.
$$
This gives
\begin{equation*}
M:=\sum_{a\in A}|B\setminus N(a)|=ab-e(A,B)\le2e(A)=O_k(n^{1+1/k}), 
\end{equation*}
recalling that $G[A]$ is $C_{2k}$-free, so $e(A)\le O(n^{1+1/k})$.
Now applying Lemma~\ref{lem:common} to $A,B$, we have 
$$
e(B)\le (\frac{k-1}{2}+\frac{kM}{a})b=O_{k,\tau}(n^{1+\frac{1}{k}}).
$$

Choose a maximum cut $A'\sqcup B'$. 
Let $a_1=|A\cap B'|$ and $b_1=|B\cap A'|$.
We assume $\frac{a_1}{a}+\frac{b_1}{b}\le1$ since otherwise we can exchange $A'$ and $B'$.
Clearly, $e(A',B')\le e(K[A,B])-L_0+e(A)+e(B)$, where $L_0=|E(K[A,B])\setminus E(K[A',B'])|$.
Also, $e(K[A,B])=ab$ and 
$$
L_0=a_1(b-b_1)+b_1(a-a_1)\ge \frac{1}{2}(ba_1+ab_1)\ge \frac{\tau n}{2}(a_1+b_1).
$$
On the other hand, $e(A',B')\ge e(A,B)=ab-M$. Hence $ab-M\le ab-L_0+e(A)+e(B)$.
Therefore,
$$
a_1+b_1\le \frac{2}{\tau n}L_0\le \frac{2}{\tau n}(M+e(A)+e(B))=O_{k,\tau}(n^{\frac{1}{k}}).
$$
So to obtain $A'\sqcup B'$ from $A\sqcup B$, we only need to move $O_{k,\tau}(n^{\frac{1}{k}})$ vertices from their partition set to the other one. 
So we still have $|A'|,|B'|\ge\tau n/2$.
The change of the vertices will cause the change of $e(A,B)$, $e(A)$ and $e(B)$ for at most $O_{k,\tau}(n^{1+\frac{1}{k}})$.
So the new cut $A'\sqcup B'$ still satisfies $|A'||B'|-e(A',B')=O_{k,\tau}(n^{1+\frac{1}{k}})$ and $e(A')+e(B')=O_{k,\tau}(n^{1+\frac{1}{k}})$.
The "moreover" part is trivial due to the maximality of $A'\sqcup B'$.
\end{proof}

The following theorem can be considered as a stability result, which guarantees a nearly dominating vertex of a $W_{2k+1}$-free graph if its $Q$-index is close to the $Q$-index of the extremal construction.

\begin{theorem}\label{thm:diffuse}
For every fixed integer $k\ge3$ and $0<\tau<1/2$, there are
constants $\xi=\xi(k,\tau)>0$ and $N=N(k,\tau)$ such that every
$n$-vertex $W_{2k+1}$-free graph $G$ with $n\ge N$ and
$\Delta(G)\le(1-\tau)n$ satisfies
$$
q(G)\le n+2(k-1)-\xi.
$$
\end{theorem}
We prove the theorem by several claims.
Throughout the proof, we consider $n$-vertex $W_{2k+1}$-free graph $G$ with $n\ge N$ and
$\Delta(G)\le(1-\tau)n$.
We will take
$$
c_1=\frac{\tau^2}{192},\qquad
\beta=\frac{c_1^2\tau}{16},
$$
and fix a positive constant $\gamma$ such that
\begin{equation}\label{eq:gamma-choice-revised}
\gamma=\frac{c_1^2\tau^3}{24k}=\min\left\{\frac{\tau}{8},\;
\frac{\tau^2}{96},\;\frac{\beta}{4k},\;
\frac{c_1^2\tau^3}{24k}\right\}.
\end{equation}
All three constants depend only on $k$ and $\tau$, and are
independent of $n$. We shall obtain a gap with
$\xi=c_1^2\tau/2$. This constant is less than $1$.
In the following lemmas we first assume $q(G)\ge n$ and choose a maximum cut $A\sqcup B$ as in Lemma~\ref{lem:cut}.
In particular,
\begin{equation}\label{eq:cut-data-revised}
|A|,|B|\ge\frac{\tau n}{2},\qquad
M:=|A||B|-e(A,B)=O_{k,\tau}(n^{1+1/k}),
\end{equation}
and $d_{B}(v_a)\ge d_{A}(v_a)$ and $d_{A}(v_b)\ge d_{B}(v_b)$ for every vertex $v_a\in A$ and $v_b\in B$.

Put $A_0=\{v\in A:|B\setminus N(v)|\le\gamma n\}$, $B_0=\{v\in B:|A\setminus N(v)|\le\gamma n\}$, $Z=V(G)\setminus(A_0\cup B_0)$, $z=|Z|$.
We also call vertices in $A_0\cup B_0$ \emph{good} vertices. Note that
$$
2M=\sum_{v\in A}|B\setminus N(v)|+\sum_{v\in B}|A\setminus N(v)|,
$$
and every vertex in $Z$ contributes more than $\gamma n$
to this sum. Hence
\begin{equation}\label{eq:z}
z\le\frac{2M}{\gamma n}=O_{k,\tau}(n^{1/k}).
\end{equation}
Since $z=o(n)$ and $k\ge3$, for sufficiently large $n$, we have
\begin{equation}\label{eq:goodsizes-revised}
|A_0|,|B_0|\ge\frac{\tau n}{3},\qquad
z^2=O_{k,\tau}(n^{2/k})=o(n).
\end{equation}
Every good vertex has at least $\tau n/2-\gamma n$
crossing neighbors, and hence has degree at least $\beta n$
for all large $n$.

\begin{claim}\label{lem:good-structure}
If $d_G(v)\ge\beta n$ for some $v\in V(G)$, then
\begin{equation*}
d_{A_0}(v)\le k-1\quad \text{if }v\in A,\quad d_{B_0}(v)\le k-1\quad \text{if }v\in B.
\end{equation*}
In particular,
\begin{equation*}
\Delta(G[A_0]), \Delta(G[B_0])\le k-1.
\end{equation*}
\end{claim}
\begin{proof}
We need to prove that every vertex $v$ with $d_G(v)\ge\beta n$
has at most $k-1$ good internal neighbors.
Indeed, since the cut is maximum, $v$ has at least
$d_G(v)/2\ge\beta n/2$ crossing neighbors.
If $v$ had $k$ good internal neighbors $s_1,\dots,s_k$,
then $v,s_1,\dots,s_k$ would have at least
$$
\frac{\beta n}{2}-k\gamma n
\ge\frac{\beta n}{4}>k
$$
common crossing neighbors. Choosing any $k$ of them and alternating them
with $s_1,\dots,s_k$ gives a $C_{2k}$ in $G[N(v)]$, which is
impossible. 
\end{proof}

Let $x$ be a nonnegative eigenvector of $Q(G)$ for $q(G)$, normalized by $\max_v x_v=1$. For $i=1,2$, set $X_i=\sum_vx_v^i$.
\begin{claim}\label{c:norm}
  $X_1\ge\tau n$ and $\tau^2n\le X_2\le n$.
\end{claim}
\begin{proof}
Let $x_u=1$. The eigenvalue equation at $u$ gives
\begin{equation*}
X_1=\sum_{v\in V(G)}x_v\ge\sum_{v\in N(u)}x_v
=q(G)-d_G(u)\ge n-\Delta(G)\ge\tau n.
\end{equation*}
By Cauchy--Schwarz inequality,
\begin{equation*}
X_2=\sum_{v\in V(G)}x_v^2\ge \frac{1}{n}\left(\sum_{v\in V(G)}x_v\right)^2\ge\tau^2n.
\end{equation*}
On the other hand, since $x_v\le1$ for $v\in V(G)$, it holds $X_2=\sum_vx_v^2\le \sum_v1=n$. The claim is done.
\end{proof}
\begin{claim}\label{c:good-positive}
For $v\in A_0\cup B_0$, $x_v\ge c_1$.
\end{claim}
\begin{proof}
Since $z=o(n)$ and $x_v\le1$ for $v\in V(G)$, by Claim~\ref{c:norm}, for large $n$,
$$
\sum_{v\in A_0\cup B_0}x_v\ge\tau n-z\ge\frac{\tau n}{2}.
$$
By symmetry we may suppose
$\sum_{v\in A_0}x_v\ge\tau n/4$.
For any $v\in B_0$, at most $\gamma n$ vertices of $A$
are nonneighbors of $v$. The eigenvalue equation gives
$$
(q(G)-d_G(v))x_v
=\sum_{w\in N(v)}x_w
\ge\sum_{w\in A_0}x_w-\gamma n
\ge\left(\frac\tau4-\gamma\right)n.
$$
Consequently, by \eqref{eq:gamma-choice-revised} and $q(G)\le2\Delta(G)<2n$ , for $v\in B_0$,
$$
x_v\ge\frac{\tau/4-\gamma}{2}\ge\frac\tau{16}.
$$
Using $|B_0|\ge\tau n/3$ yields
$\sum_{v\in B_0}x_v\ge\tau^2n/48$.
Applying the same eigenvalue estimate to each $v\in A_0$ gives
$$
x_v\ge\frac{\tau^2/48-\gamma}{2}
\ge\frac{\tau^2}{192}.
$$
Since $c_1=\tau^2/192$, we have proved that for $v\in A_0\cup B_0$,
$x_v\ge c_1$.
\end{proof}

Let $I_0=\sum_{uv\in E(G[A_0])\cup E(G[B_0])}(x_u+x_v)^2$.
\begin{claim}\label{lem:good-energy}
$
q(G)\le n+\frac{I_0}{X_2}+o_{k,\tau}(1).
$
\end{claim}
\begin{proof}
Using the identity
$$
x^{\mathsf T}Q(J)x=\sum_{uv\in E(J)}(x_u+x_v)^2
$$
for any graph $J$, we have
\begin{equation}\label{eq:energy-identity-revised}
\begin{aligned}
x^{\mathsf T}Q(G)x
={}&x^{\mathsf T}Q(K[A,B])x+\sum_{uv\in E(G[A])\cup E(G[B])}(x_u+x_v)^2\\
&-\sum_{\substack{u\in A,\ v\in B\\uv\notin E(G)}}(x_u+x_v)^2.
\end{aligned}
\end{equation}
Note that $x^{\mathsf T}Q(K[A,B])x\le n X_2$ as $Q(K[A,B])$ has largest eigenvalue $n$.

Let $Z_{\rm low}=\{v\in Z:d_G(v)<\beta n\},~z_0=|Z_{\rm low}|.$
We now consider the terms in \eqref{eq:energy-identity-revised}.

Since $(x_u+x_v)^2\le4$,
all internal edges with one endpoint in $Z_{\rm low}$ and the
other in $A_0\cup B_0$ contribute at most $4\beta nz_0$.
For any $v\in Z_{\rm low}$, at least
$$
\frac{\tau n}{3}-d_G(v)
\ge\left(\frac\tau3-\beta\right)n
\ge\frac{\tau n}{4}
$$
good vertices are crossing nonneighbors of $v$.
Each of these missing edges contributes at least $c_1^2$
to the negative sum by Claim~\ref{c:good-positive}.
Thus the total negative contribution of all crossing nonedges with one endpoint in $Z_{\rm low}$ and the other in $A_0\cup B_0$ is at least
$c_1^2\tau nz_0/4=4\beta nz_0$. Together with the positive contribution of the internal edges with one endpoint in $Z_{\rm low}$ and the
other in $A_0\cup B_0$ we counted before, we now have the total contribution at most $4\beta nz_0-4\beta nz_0=0$.

The remaining internal edges not counted by $I_0$ but by \eqref{eq:energy-identity-revised} are of two types.
Those with one endpoint in $A_0\cup B_0$ and the other in $Z\setminus Z_{\rm low}$ contribute at most
$4(k-1)z$ by Claim~\ref{lem:good-structure}.
Those with both endpoints in $Z$ contribute at most
$4\binom z2\le2z^2$.
All other terms from missing crossing edges are nonpositive.
Consequently \eqref{eq:energy-identity-revised} yields
$$
q(G) X_2\le n X_2+I_0+4(k-1)z+2z^2.
$$
Using \eqref{eq:z}, \eqref{eq:goodsizes-revised} and
Claim~\ref{c:norm}, we conclude
\begin{equation}\label{eq:energy}
q(G)\le n+\frac{I_0}{X_2}+o(1).
\end{equation}
Indeed the omitted error is bounded by
$$
\frac{4(k-1)z+2z^2}{X_2}
=O_{k,\tau}(n^{1/k-1}+n^{2/k-1})=o_{k,\tau}(1).
$$
\end{proof}

\medskip
\begin{proof}[Proof of Theorem~\ref{thm:diffuse}]
If $q(G)<n$, the desired inequality is immediate. Assume $q(G)\ge n$.
Suppose first that $\Delta(G[A_0]), \Delta(G[B_0])\le k-2$.
Since $(x_u+x_v)^2\le2(x_u^2+x_v^2)$,
$$
I_0\le2(k-2)\sum_{v\in A_0\cup B_0}x_v^2
\le2(k-2)X_2.
$$
By Claim~\ref{lem:good-energy},
$$
q(G)\le n+2(k-2)+o(1)=n+2(k-1)-2+o(1),
$$
which gives the required fixed gap for large $n$.

Otherwise, by symmetry and Claim~\ref{lem:good-structure},
there is a vertex $u\in A_0$ with exactly $k-1$ neighbors
in $A_0$. Put
$$
T=B_0\cap N(u)\cap\bigcap_{v\in N_{A_0}(u)}N(v).
$$
Both $u$ and all vertices in $N_{A_0}(u)$ are good, so
\begin{equation}\label{eq:T-complement-revised}
|B_0\setminus T|\le k\gamma n.
\end{equation}
Since $|B_0|\ge\tau n/3$ and
$k\gamma\le\beta/4<\tau/12$, we have
$|T|\ge\tau n/4\ge k+1$ for all large $n$.
The sets $N_{A_0}(u),T$ and the vertex $u$ therefore satisfy the hypotheses
of Lemma~\ref{lem:star} (1). 
It follows that $e(T)\le1$.
By Claim~\ref{lem:good-structure} and
\eqref{eq:T-complement-revised},
\begin{equation*}
e(B_0)\le e(T)+(k-1)|B_0\setminus T|
\le1+k(k-1)\gamma n.
\end{equation*}
Using the degree bound in $G[A_0]$ and
$(x_u+x_v)^2\le2(x_u^2+x_v^2)$, we have
\begin{equation}\label{eq:i0}
\begin{aligned}
I_0
&\le2(k-1)\sum_{v\in A_0}x_v^2+4e(B_0)\\
&\le2(k-1)X_2
 -2(k-1)\sum_{v\in B_0}x_v^2
 +4+4k(k-1)\gamma n.
\end{aligned}
\end{equation}
On the other hand, Claim~\ref{c:good-positive} and
\eqref{eq:goodsizes-revised} imply
$$
\sum_{v\in B_0}x_v^2\ge c_1^2|B_0|
\ge\frac{c_1^2\tau n}{3}.
$$
Combining this with \eqref{eq:i0}, Claim~\ref{c:norm} and Claim~\ref{lem:good-energy}, we obtain
\begin{align*}
q(G)
&\le n+2(k-1)-\frac{2(k-1)c_1^2\tau}{3}
+\frac{4k(k-1)\gamma}{\tau^2}+o(1)\\
&= n+2(k-1)-\frac{(k-1)c_1^2\tau}{2}+o(1),
\end{align*}
where the last equality uses
$\gamma=c_1^2\tau^3/(24k)$.
Since $k,\tau$ and all constants have been fixed before
$n\to\infty$, by increasing $N(k,\tau)$ if necessary we
obtain in either case
$$
q(G)\le n+2(k-1)-\xi,
\qquad\xi=\frac{c_1^2\tau}{2}>0.
$$
This proves the theorem.
\end{proof}

\section{Proof of Theorem~\ref{thm:main}}
\begin{proof}[Proof of Theorem~\ref{thm:main}]
Assume that Theorem~\ref{thm:main} is false. Since $S_{n,k}^{+}$ is $W_{2k+1}$-free, by Lemma~\ref{lem:benchmark}, there are arbitrarily large integers $n$ and $n$-vertex $W_{2k+1}$-free graphs $G\not\cong S_{n,k}^{+}$ such that
$$
n+2(k-1)-\frac{2k(k-1)}n< q(S_{n,k}^{+})\le q(G).
$$

We first show that
\begin{equation}\label{eq:asymptotic-universal-revised}
\frac{n-1-\Delta(G)}{n}\longrightarrow0.
\end{equation}
Indeed, fix any $0<\tau<1/2$.
By Theorem~\ref{thm:diffuse}, for all sufficiently large
$n$, if $\Delta(G)\le(1-\tau)n$, $q(G)\le n+2(k-1)-\xi(k,\tau)$, where $\xi(k,\tau)>0$ does not depend on $n$.
This contradicts the preceding lower bound on $q(G)$ when $n$ is sufficiently large. Thus  $\Delta(G)>(1-\tau)n$ for every fixed $\tau>0$ small enough. Since $\Delta(G)\le n-1$, this proves
\eqref{eq:asymptotic-universal-revised}.

Now fix $\eps=1/4$ and use the notations before:
$$
\alpha=\frac{\eps}{2k},\quad
L=L_{G}(\eps n),\quad
U=L_{G}(n-1-\alpha n),\quad
\delta(v)=n-1-d_{G}(v).
$$
Choose a vertex $u$ of maximum degree in $G$.
Equation~\eqref{eq:asymptotic-universal-revised} gives
$$
\delta(u)=n-1-\Delta(G)=o(n).
$$
Since $\alpha>0$ is fixed, for all sufficiently large $n$
we have $\delta(u)\le\alpha n$, that is, $u\in U$.
We may also require $n\ge N_0(k,\eps)$.
All hypotheses of Lemma~\ref{lem:nearu} (3) are therefore satisfied,
which gives
$$
d_L(u)\le D_0=\left\lceil(k-1)(4/\eps)^k\right\rceil.
$$
Every vertex of $L$ is either $u$, a neighbor of $u$ in $L$,
or a nonneighbor of $u$. The number of nonneighbors of $u$
is at most $\delta(u)$. Therefore,
$$
|L|\le1+d_L(u)+\delta(u)\le1+D_0+o(n)=o(n).
$$
We can now apply Theorem~\ref{thm:sublinear} to $G$ with this fixed $\eps$. It follows that $G\cong S_{n,k}^{+}$ for all sufficiently large $n$, contradicting our choice of $G$.
Thus $q(G)\le q(S_{n,k}^{+})$ for every sufficiently large $n$-vertex $W_{2k+1}$-free graph $G$, and the equality case follows from Theorem~\ref{thm:sublinear} by the same argument.
\end{proof}

\section{Concluding remarks}
In this paper, we show that $q(G)\le q(S_{n,k}^+)$, with equality if and only if $G\cong S_{n,k}^+$, for any $n$-vertex $W_{2k+1}$-free graph $G$ when $n$ is large enough and $k\ge 3$.
However, for $k=2$, the problem is still open. Our proof does not work because Lemma~\ref{lem:star} fails for $k=2$ and the extremal construction is totally different. 

For integers $a,b\ge 2$, let 
$H_{a,b}$ be the graph obtained from $K_{a,b}$ by adding a maximum matching to each of its partition sets. To be more specific, 
$$H_{a,b}:= 
\left(\left\lfloor\frac{a}{2}\right\rfloor K_2\cup \left(a-2\left\lfloor\frac{a}{2}\right\rfloor\right)K_1\right)
\vee 
\left(\left\lfloor\frac{b}{2}\right\rfloor K_2\cup \left(b-2\left\lfloor\frac{b}{2}\right\rfloor\right)K_1\right).$$
If $a$ and $b$ are both even numbers, $H_{a,b}=\frac{a}{2}K_2\vee\frac{b}{2}K_2$. One can check that $q(H_{a,b})=a+b+2$ in this case.
We give the following conjecture.
\begin{conjecture}
There is a constant $N$ such that, for any $n\ge N$, every $n$-vertex $W_{5}$-free graph $G$ satisfies $q(G)\le q(H_{2,n-2})$. 
For odd $n$, the equality holds if and only if $G\cong H_{2,n-2}$.
For even $n$, the equality holds if and only if $G\cong H_{a,b}$ for some positive even numbers $a$ and $b$ with $a+b=n$.
\end{conjecture}

\subsection*{Declaration of generative AI use}
During the research and preparation of this manuscript, the authors used ChatGPT (OpenAI) to explore possible proof strategies, check mathematical arguments and calculations, and improve the organization and language of the manuscript. The authors critically reviewed and revised the AI-assisted material and take full responsibility for the mathematical correctness and content of the paper.

\subsection*{Acknowledgements}
This work was supported by National Natural Science Foundation of China (Grant No. 12401455, 12571369) and the Stability Support Program for Basic Research (Grant No. WZJC202601001).


\begin{thebibliography}{99}

\bibitem{BondySimonovits1974}
J.A. Bondy, M. Simonovits,
Cycles of even length in graphs,
J. Combin. Theory Ser. B 16 (1974) 97--105.

\bibitem{ChenMa2023}
D. Chen, X. Ma,
The maximum signless Laplacian spectral radius
of graphs with forbidden subgraphs,
Filomat 37 (2023) 8319--8330.

\bibitem{ChenGaoLi2026}
H. Chen, J. Gao, S. Li,
The edge spectral extremal problem for odd wheels in nonzero residue classes,
arXiv:2609.13666 (2026).

\bibitem{ChenJinZhang2025}
M. Chen, Y. Jin, P. Zhang,
The signless Laplacian spectral radius of book-free graphs,
Discrete Math. 348 (2025) 114448.

\bibitem{ChenJinZhang2026}
M. Chen, Y. Jin, P. Zhang,
The signless Laplacian spectral Tur\'an problems for disjoint triangles,
arXiv:2610.30859 (2026).

\bibitem{ChenLiuZhang2020}
M. Chen, A. Liu, X. Zhang,
The signless Laplacian spectral radius of graphs with forbidding linear forests,
Linear Algebra Appl. 591 (2020) 25--43.

\bibitem{ChenLiuZhang2023}
M. Chen, A. Liu, X. Zhang,
On the $A_\alpha$-spectral radius of graphs without
linear forests,
Appl. Math. Comput. 450 (2023) 128005.

\bibitem{ChenLiuZhang2024}
M. Chen, A. Liu, X. Zhang,
The signless Laplacian spectral radius of graphs without intersecting odd cycles,
Electron. J. Linear Algebra 40 (2024) 370--381.

\bibitem{CioabaDesaiTait2022}
S.M. Cioab\u{a}, D.N. Desai, M. Tait,
The spectral radius of graphs with no odd wheels,
European J. Combin. 99 (2022) 103420.

\bibitem{FengGu2009}
L. Feng, G. Yu, On three conjectures involving the signless Laplacian spectral radius of graphs, Publ. Inst. Math. (Beograd) 85(99) (2009), 35–38.

\bibitem{FreitasNikiforovPatuzzi2013}
M.A.A. de Freitas, V. Nikiforov, L. Patuzzi,
Maxima of the $Q$-index: Forbidden 4-cycle and 5-cycle,
Electron. J. Linear Algebra 26 (2013) 905--916.

\bibitem{FreitasNikiforovPatuzzi2016}
M.A.A. de Freitas, V. Nikiforov, L. Patuzzi,
Maxima of the $Q$-index: Graphs with no $K_{s,t}$,
Linear Algebra Appl. 496 (2016) 381--391.

\bibitem{HeJinZhang2013}
B. He, Y. Jin, X. Zhang, Sharp bounds for the signless Laplacian spectral radius in terms of clique number, Linear Algebra Appl. 438 (2013) 3851--3861.

\bibitem{Merris1998}
R. Merris,
A note on Laplacian graph eigenvalues,
Linear Algebra Appl. 285 (1998) 33--35.

\bibitem{NikiforovYuan2014}
V. Nikiforov, X. Yuan,
Maxima of the $Q$-index: Graphs without long paths,
Electron. J. Linear Algebra 27 (2014) 504--514.

\bibitem{NikiforovYuan2015}
V. Nikiforov, X. Yuan,
Maxima of the $Q$-index: Forbidden even cycles,
Linear Algebra Appl. 471 (2015) 636--653.

\bibitem{WangZhai2023}
B. Wang, M. Zhai,
Maxima of the $Q$-index: Forbidden a fan,
Discrete Math. 346 (2023) 113264.

\bibitem{Yuan2014}
X. Yuan,
Maxima of the $Q$-index: Forbidden odd cycles,
Linear Algebra Appl. 458 (2014) 207--216.

\bibitem{YuanZhaiZhang2026}
W. Yuan, M. Zhai, Y. Zhang,
Maximal $Q$-index: Graphs without long wheels,
Linear Algebra Appl. 749 (2026) 1--14.

\bibitem{Zhang2025}
W. Zhang,
More results on the spectral radius of graphs with no odd wheels,
Discrete Math. 348 (2025) 114550.

\bibitem{ZhangWang2024}
Y. Zhang, L. Wang,
The signless Laplacian spectral radius of $2K_3$-free graphs,
Discrete Math. 347 (2024) 114075.

\bibitem{ZhaoHuangGuo2021}
Y. Zhao, X. Huang, H. Guo,
The signless Laplacian spectral radius of graphs with no intersecting triangles,
Linear Algebra Appl. 618 (2021) 12--21.

\bibitem{ZhaoHuangLin2021}
Y. Zhao, X. Huang, H. Lin,
The maximum spectral radius of wheel-free graphs,
Discrete Math. 344 (2021) 112341.

\bibitem{ZhaoYouZeng2026}
X. Zhao, L. You, J. Zeng,
The signless Laplacian spectral radius of graphs without disjoint cliques,
arXiv:2606.30272 (2026).

\bibitem{ZhengLiLi2026}
J. Zheng, Y. Li, H. Li,
The signless Laplacian spectral Tur\'an problems for  color-critical graphs,
Linear Algebra Appl. 730 (2026) 546--565.

\bibitem{ZhengLiSu2026}
J. Zheng, H. Li, L. Su,
A signless Laplacian spectral Erd\H{o}s--Stone--Simonovits theorem,
Discrete Math. 349 (2026) 114665.

\end{thebibliography}
\end{document}